\RequirePackage[OT1]{fontenc}

\documentclass[letterpaper, 10 pt, conference]{ieeeconf}
\IEEEoverridecommandlockouts

\usepackage{cite}
\usepackage{amsmath,amssymb,amsfonts}
\usepackage{algorithmic}
\usepackage{graphicx}
\usepackage{textcomp}
\usepackage{xcolor}
\usepackage{mathtools}
\usepackage{bm}
\usepackage{optidef}
\usepackage{pgf}
\usepackage{float}
\usepackage{booktabs}
\usepackage{multirow}
\usepackage{tabularx}

\newcommand{\rr}{\mathbb{R}}
\newcommand{\nn}{\mathbb{N}}
\newcommand{\Ll}{L}
\newcommand{\Uu}{\mathcal{U}}

\newcommand{\Xx}{\mathcal{X}}
\newcommand{\Pp}{\mathcal{P}}
\newcommand{\nuu}{n_{u}}

\newcommand{\nq}{n_{q}}

\newcommand{\np}{n_{p}}
\newcommand{\nint}{n_{\mathrm{int}}}
\newcommand{\nc}{n_{g}}
\newcommand{\nr}{n_{r}}

\newcommand{\dt}{\mathrm{d}t}

\def\BibTeX{{\rm B\kern-.05em{\sc i\kern-.025em b}\kern-.08em
		T\kern-.1667em\lower.7ex\hbox{E}\kern-.125emX}}

\newcolumntype{L}[1]{>{\raggedright\arraybackslash}p{#1}}
\newcolumntype{C}[1]{>{\centering\arraybackslash}p{#1}}
\newcolumntype{R}[1]{>{\raggedleft\arraybackslash}p{#1}}

\newcommand\copyrighttext{%
	\footnotesize \textcopyright 2022 IEEE. Personal use of this material is permitted. Permission from IEEE must be obtained for all other uses, in any current or future media, including reprinting/republishing this material for advertising or promotional purposes, creating new collective works, for resale or redistribution to servers or lists, or reuse of any copyrighted component of this work in other works.}
\usepackage{tikz}
\newcommand\copyrightnotice{%
	\begin{tikzpicture}[remember picture,overlay]
		\node[anchor=south,yshift=5pt] at (current page.south) {\fbox{\parbox{\dimexpr\textwidth-\fboxsep-\fboxrule\relax}{\copyrighttext}}};
	\end{tikzpicture}%
	\thispagestyle{empty}
}

\begin{document}
	
	\title{\LARGE \bf Direct Collocation for Numerical Optimal Control of Second-Order ODE  \\
		\thanks{This research was supported by the DFG via Research Unit FOR 2401 and project 424107692 and by the EU via ELO-X 953348.
		}
		\thanks{L\'eo Simpson is with the Research and Development team, Tool-Temp AG, Switzerland,
			\texttt{leo.simpson@tool-temp.ch. }}
		\thanks{Armin Nurkanovi\'c is with the Department of Microsystems Engineering (IMTEK), University of Freiburg, Germany, 
			\texttt{armin.nurkanovic@imtek.uni-freiburg.de. }}
		\thanks{Moritz Diehl is with the
			Department of Microsystems Engineering (IMTEK) and Department of Mathematics, University of Freiburg, Germany,  \texttt{armin.nurkanovic,moritz.diehl@imtek.uni-freiburg.de}
		}
	}
	
	\author{L\'eo Simpson, Armin Nurkanovi\'c, Moritz Diehl}
	
	\maketitle
	\copyrightnotice
	
	\begin{abstract}
		Mechanical systems are usually modeled by second-order Ordinary Differential Equations (ODE) which take the form $\ddot{q} = f(t, q, \dot{q})$.
		While simulation methods tailored to these equations have been studied, using them in direct optimal control methods is rare.
		Indeed, the standard approach is to perform a state augmentation, adding the velocities to the state.
		The main drawback of this approach is that the number of decision variables is doubled, which could harm the performance of the resulting optimization problem.
		In this paper, we present an approach tailored to second-order ODE.
		We compare it with the standard one, both on theoretical aspects and in a numerical example.
		Notably, we show that the tailored formulation is likely to improve the performance of a direct collocation method, for solving optimal control problems with second-order ODE of the more restrictive form $\ddot{q} = f(t, q)$.
		
	\end{abstract}
	
	\section{Introduction}
	
	Optimal planning of dynamic trajectories in robotics requires an efficient way of solving finite-horizon optimal control problems.
	In particular, online planning in the context of nonlinear model predictive control requires repeated and fast solutions to these problems \cite{magni2009nonlinear}.
One way of solving such problems is direct optimal control, which starts with a discretization of the continuous problem before solving it.
Among the direct methods, direct collocation methods are widely used, as they generate sparse Nonlinear Programming (NLP) problems that can be solved efficiently with dedicated solvers \cite{rawlings2017model}.
	
	This paper is concerned with Optimal Control Problems (OCP) with second-order Ordinary Differential Equations (ODE), which read as: 
	\begin{align}\label{2ndorderODE-control}
		\ddot{q}(t) & = f(q(t), \dot{q}(t), u(t), p), \quad t \in [0, T],
	\end{align}
	where $ T \in \rr$,  $q : [0, T] \to \rr^{\nq}$, $\dot{q}(\cdot)$ and $\ddot{q}(\cdot)$ are the first and second time derivative of $q(\cdot)$, $u : [0, T] \to \rr^{\nuu}$, and $f : \rr^{\nq} \times \rr^{\nq} \times \rr^{\nuu} \times \rr^{\np} \to \rr^{\nq}$ is a Lipschitz function.
	This structure of ODE is very common.
It can model dynamical systems subject to classical mechanical laws \cite{lurie2002analytical}.
	For example, it can model the equations of motion for walking robots \cite{winkler2018gait}.
	Simulation methods tailored to this structure have been studied in \cite{kramarz1980stability, van1991stability, d2009two, martucci2005general}.
	In this paper, we investigate how to adapt  numerical optimal control methods to the presented special structures to increase computational efficiency.
	In recent work \cite{moreno2022collocation}, two direct collocation methods for optimal control are adapted to second-order ODE: the trapezoidal and the Hermite-Simpson collocation methods.
	These methods only imply low-order convergence compared to Legendre-Gauss-Radau collocation methods \cite{rawlings2017model}.
	For these two methods, improvements compared to the standard methods have been shown for two numerical examples \cite{moreno2022collocation}.

	\paragraph*{Contributions}
	In this paper, we present and compare two approaches to systematically apply collocation methods to the second-order structure \eqref{2ndorderODE-control} in optimal control settings.
	The first approach is standard for direct optimal control, but does not exploit the special structure while the second approach does.
	More precisely, in the first approach, the ODE \eqref{2ndorderODE-control} is considered as a first-order ODE where the state comprises positions and velocities, while in the second, one directly considers the second-order structure.
	We compare the two approaches theoretically and numerically.
	The authors \cite{moreno2022collocation} exploit the aforementioned structure in two special cases.
This paper extends these to any collocation method, and in particular, to the higher-order ones.
	To formulate these new direct collocation methods, we use a special polynomial basis that permits us to formulate new direct collocation methods with more sparsity in the equations.
	We show that for the more restrictive case $\ddot{q} = f(q, u)$ the NLP resulting from our new method is less complex (lower dimension and similar sparsity pattern) than for the standard approach and yet the underlying integration method is as accurate, as shown on a numerical example.

	\paragraph*{Outline}
	In Section \ref{section:collocation} we define two classes of methods to apply collocation methods to second-order systems which we call Standard Collocation (SC) and Position-based Collocation (PC).
	In Section \ref{section:direct-collocation} we define the class of OCP we tackle and present direct collocation to solve these numerically, using the integration methods defined before.
	This is followed by Section \ref{section:theoretical}, where an advantage of the new method is highlighted.
	It is finally put into practice in Section \ref{section:numerical} on numerical examples.
	
	\paragraph*{Notations}
	In this paper, we denote $\left(a_1, ..., a_n \right) \coloneqq \left[a_1^{\top}, ..., a_n^{\top}\right]^{\top}$ the concatenated column vector when $a_1, ..., a_n$ are column vectors.

	\section{Collocation methods for second-order systems}\label{section:collocation}
	In this section, we are not concerned with optimal control yet, but only with integration methods, hence we will consider a second-order ODE in the following form:
	\begin{align}\label{2ndorderODE}
		\ddot{q}(t) & = f(t, q(t), \dot{q}(t)), \quad t \in [0, T],
	\end{align}
	where $f : \rr \times \rr^{\nq} \times \rr^{\nq} \to \rr^{\nq}$ is a Lipschitz function.
	
	Furthermore, we make use of the change of variable $x \coloneqq \left(q, \dot{q}\right)$
	and define the corresponding function
	$\bar{f}\left( t, x\right) \coloneqq \left( \dot{q}, \;  f(t, q, \dot{q}) \right)$.
	With this change of variable, the second-order ODE \eqref{2ndorderODE} becomes the first-order ODE:
	\begin{align}\label{ODE}
		\dot{x}(t) & = \bar{f}(t, x(t)), \quad t \in [0, T].
	\end{align}
	
	We use an equidistant grid
	defined by the grid size $N$ and the grid-points $t_{k} \coloneqq kh$ for $k \in \mathcal{N} \coloneqq \{0, ..., N\}$, with $h \coloneqq \frac{T}{N}$ the step-size.

	\subsection{Collocation methods for second-order ODE}\label{subsection:methods}

	In this section, we are concerned with the second-order Initial Value Problem (IVP) given by the ODE \eqref{ODE} and the initial condition $x(0) = x_0$ where  $ x_0 \coloneqq (q_0, v_0)$ is a vector in $\rr^{2\nq}$.
	We also fix $d\in \nn$ and $\tau_0, ..., \tau_d \in (0, 1]$ some distinct numbers which we call \textit{collocation order} and \textit{collocation points} respectively.
	We also define the set $\mathcal{D} \coloneqq \{1, ..., d\}$.
	To derive the formulas with more simplicity we also define $\tau_0 = 0$.
	Finally, we define the corresponding points $t_{k, i} \coloneqq t_k + h \tau_{i} \in [0, T]$
	
	Collocation methods are systems of equations that come from the discretization of this IVP.
	They are originally used for simulation.
	They fall into the class of Implicit Runge-Kutta methods \cite{rawlings2017model}, but are formulated differently: they consist in writing the conditions for a piece-wise polynomial function to satisfy the differential equations on the collocation points.
	In \cite{van1991stability}, the author distinguishes two types of collocation methods for simulation of systems in the form \eqref{2ndorderODE}: \textit{indirect} and \textit{direct} simulation methods.
	In the first one, the collocation method is applied to the IVP in standard form \eqref{ODE},
	i.e. after the change of variable defined above, while in the second kind it is applied directly to the second-order ODE \eqref{2ndorderODE}.
	To avoid conflicts of terminology with the name \textit{direct collocation methods} used in the numerical optimal control community, for example in \cite{rawlings2017model}, 
	in this paper, we use different names: the first methods ("indirect" in \eqref{2ndorderODE}) will be called Standard Collocation (SC)  because it is the usual approach to direct optimal control of second-order ODE and the second ("direct" in \eqref{2ndorderODE}) Position-based Collocation (PC) .
	
	The solution $x(\cdot)$ is approximated by a polynomial on each interval $[t_k, t_{k+1}]$.
	The variable $x_k \coloneqq (q_k, v_k) \in \rr^{2\nq}$ approximates $x(t_k) $ and $z_k \in \rr^{\nint}$ are intermediate variables which contain further information about the solution within the interval $[t_k, t_{k+1}]$, with a certain $\nint \in \nn$ depending on the chosen collocation method.
	In standard collocation methods, the approximation is of the form $x(t_k + h \tau) \simeq \phi(x_k, z_k, \tau)$ for $\tau \in [0, 1]$ where $\phi(x_k, z_k, \cdot)$ is a polynomial with $2\nq$ components of degree $d$ or less, and is parameterized by the intermediate variables $z_k \in \rr^{2\nq \cdot \, d}$ through a well-chosen polynomial basis such that $\phi_k(x_k, z_k, 0) = x_k$.
	We denote by $\dot{\phi}(x_k, z_k, \tau)$ the derivative of this polynomial with respect to $\tau$.
	The standard collocation method is given by the following recursive equations for $k \in \mathcal{N}$:
	\begin{align}\label{eq:standard}
		\begin{split}
			x_{k+1} &= \phi(x_k, z_k, 1),  \\
			\dot{\phi}(x_k, z_k, \tau_i) &= h\bar{f} \left(t_{k,i}, \; \phi(x_k, z_k, \tau_i) \; \right) \hspace{0.5cm}\mbox{for all} \;i \in \mathcal{D}.
		\end{split}		
	\end{align}

	In PC methods, the approximation is of the form:
	\begin{align}
		x(t_k + h \tau) \simeq \left( \psi(x_k, z_k, \tau), \; h^{-1}\dot{\psi}(x_k, z_k, \tau) \right),
	\end{align}
	for $\tau \in [0, 1]$.
Furthermore, $\psi(x_k, z_k, \cdot)$ is a polynomial with only $\nq$ components of degree $d+1$ or less, and is parameterized by the intermediate variables $z_k \in \rr^{\nq \cdot \, d}$ through a well-chosen polynomial basis such that $x_k = \left( \psi(x_k, z_k, 0), h^{-1}\dot{\psi}(x_k, z_k, 0) \right)$.
	
	The PC method is given by the following recursive equations for $k \in \mathcal{N}$:
	\begin{align}\label{eq:acceleration}
		\begin{split}
			x_{k+1} &= \left( \psi(x_k, z_k, 1), \; h^{-1}\dot{\psi}(x_k, z_k, 1) \right),  \\
			\ddot{\psi}(x_k, z_k, \tau_i) &= h^2 f\left(t_{k,i}, \; \psi(x_k, z_k, \tau_i), \; h^{-1}\dot{\psi}(x_k, z_k, \tau_i)  \; \right) \\ 
			&\hspace{5cm}  \mbox{for all} \;i \in \mathcal{D}.
		\end{split}
	\end{align}
	
	Note that both of these methods can be summarized in the following form:
	\begin{align}\label{eq:collocation-stacked}
		\begin{split}
			x_{k+1} - F(x_k, z_k) &= 0,  \\
			G_k(x_k, z_k) &= 0
		\end{split}	
	\end{align}

	It has been shown in \cite{van1991stability} that the latter method is equivalent to a Runge-Kutta-Nystr{\"o}m integration scheme and that both methods have the same order of convergence of the global error with respect to the grid size.
	Note that there is the same number of evaluations of $f(\cdot)$ in both methods, but the number of intermediate variables for standard collocation methods is $\nint = 2 d \nq$ while for PC methods, we have only $\nint = d \nq$.
	The reduction of the number of equations and degrees of freedom in \eqref{eq:acceleration} is further discussed in Subsection \ref{subsection:sparsity}.

	\subsection{Interpretation of the two methods}
	
	In this section, we discuss the difference between the methods.
	For this purpose, we start by 
	noting $\tilde{q}(t), \tilde{v}(t)$ the piece-wise polynomial functions for position and velocities that a given collocation method yields.
	
	Let us remark that for both methods, the equation $\dot{\tilde{v}}(t) = f(t, \tilde{q}(t), \tilde{v}(t)) $ holds only for $t= t_{k,i}$ with $(k, i) \in \mathcal{N} \times \mathcal{D}$.
	On the other hand, for standard collocation methods, the equation $\dot{\tilde{q}}(t) = \tilde{v}(t)$ also only holds for $t$ being a collocation point, while it does for every $t$ in PC methods.
	This difference comes from the fact that in standard methods, the equation $\dot{\tilde{q}}(t) = \tilde{v}(t)$ is imposed at $d$ points which is not sufficient to impose total equality of the polynomials, as $\tilde{v}(t)$ can be of degree $d$, which makes $d+1$ degrees of freedom.
	In \cite{moreno2022collocation}, a discussion is brought about the consequences of this fact.
It is highlighted that in standard methods, because the equation $\dot{\tilde{q}}(t) = \tilde{v}(t)$  is not satisfied everywhere, the equation $\ddot{\tilde{q}}(t) = \dot{\tilde{v}}(t)$ is usually not satisfied at any point.
As a consequence, the equation of the dynamic $\ddot{\tilde{q}}(t) = f( t, \tilde{q}(t), \dot{\tilde{q}}) $ is in general not satisfied.
It is however the case, at least on the collocation points, in PC methods.
This phenomenon is illustrated for two numerical examples in \cite{moreno2022collocation}.
	
	While this difference seems to imply a higher accuracy in PC methods, the order of convergence of the global error is the same, because the approximation error on $\dot{q}$ dominates the error on $q$, as we will see in Subsection \ref{subsection:accuracy}.

	\subsection{Adapted choice of polynomial basis}
	
	In this section, we define an adapted polynomial basis in which the PC method \eqref{eq:acceleration} can be expressed explicitly.
	Indeed, while the Lagrange polynomial basis is the obvious choice for SC methods, to apply the PC, we need to choose a parameterization of the polynomials $\psi(x_k, z_k, \cdot)$ introduced earlier, which should satisfy the equation $\left( \psi(x_k, z_k, 0), \; h^{-1}\dot{\psi}(x_k, z_k, 0) \right) = x_k$.
	
	The choice of parameterization, which defines the intermediate variables, does not affect the accuracy of the method,
	but the sparsity of the equations, which are viewed as constraints in direct optimal control.
Thus, the equations of the method have to be as simple as possible.
	To make this parameterization efficient, we define an adapted polynomial basis as follows.
	
	We define $p^{*}_0, ..., p^{*}_{d},  p^{*}_{\mathrm{v}}$ as the unique polynomials of degree less or equal to $d+1$ that satisfy the following relations:
	\begin{align}\label{eq:lagrange}
		p^{*}_i(\tau_{i'}) &= 
		\begin{cases}
			1& \text{if } i = i', \\
			0             & \text{otherwise,}
		\end{cases}\; &&\mbox{for all} \; (i, i') \in \{0, ..., d\}^2, \nonumber \\
		\dot{p}^{*}_i(0) &= 0, \; &&\mbox{for all} \; i \in \{0, ..., d\}, \\
		p^{*}_{v}(\tau_i) &= 0, \; &&\mbox{for all} \; i \in \{0, ..., d\}, \nonumber\\
		\dot{p}^{*}_{v}(0) &= 1. \nonumber
	\end{align}
	
	The polynomials $p^{*}_0, ..., p^{*}_{d},  p^{*}_{v}$ define a basis of $\rr_{d+1}[X]$ which we call the semi-Hermite polynomial basis.
	Using the polynomial $r(\tau) \coloneqq  \prod\limits_{j=1}^d  \frac{\tau_j - \tau }{\tau_j}$, these polynomials can be written explicitly:
	\begin{align*}
		p^{*}_{0}(\tau) &\coloneqq  (1 - r'(0)\tau)r(\tau), \\
		p^{*}_{v}(\tau) &\coloneqq \tau r(\tau),\\
		p^{*}_j(\tau) &\coloneqq   \frac{\tau^2}{\tau_j^2}\prod\limits_{i\in \{1, ..., d\} \setminus \{i\}}  \frac{ \tau - \tau_{i}}{\tau_{i} - \tau_{i}}, \; \mbox{for all} \; j \in \{1, ..., d\}.
	\end{align*}

\begin{figure}
		\vspace{0.5cm}
		\begin{center}
			\scalebox{0.8}{
				\input{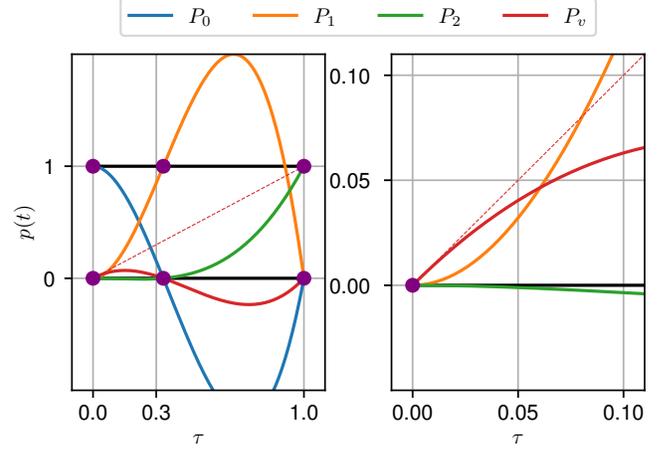}
			}
		\end{center}
		\vspace{-0.5cm}
		\caption{Semi-Hermite polynomial basis for $d=2$, and the Legendre-Radau collocation roots $(\tau_1, \tau_2) \coloneqq (0.33, 1)$.
		The plot on the right is a zoom-in of the left plot to have a better vision on $\dot{p}(0)$.}
		\label{fig:lagrange-hermite}
	\end{figure}

	Indeed, 
		direct computations show that these polynomials verify the equalities \eqref{eq:lagrange} using $r(0) = 1$ and  $r(\tau_i) = 0$ for $i\in \{1, ..., d\}$.
	These polynomials are illustrated in Figure \ref{fig:lagrange-hermite}.
	
	\subsection{Explicit formulations}
	
	For standard collocation methods, we can express the system \eqref{eq:standard} explicitly as follows:
	\begin{align}\label{eq:standard-lagrange}
		\begin{split}
			x_{k+1} &= x_k p_0(1) + \sum_{j \in \mathcal{D}} z_{k,j} p_j(1),  \\
			x_k \dot{p}_0(\tau_i) + \sum_{j \in \mathcal{D}} z_{k,j} \dot{p}_j(\tau_i) &= h \bar{f}(t_k+h\tau_{i}, z_{k,j}) \hspace{0.1cm}  \mbox{for all} \;i \in \mathcal{D}.
		\end{split}	
	\end{align}
	where $p_0, ..., p_d$ is the Lagrange polynomial basis with respect to the points $\tau_0, ..., \tau_d$.

	For PC methods, to facilitate the explicit formulation of the system \eqref{eq:acceleration}, we will start with the restricted case where $f$ does not take $\dot{q}$ as input but only $q$ and $t$.
The method can still be extended to the general case, even though the equations would be more complex, as briefly mentioned in Subsection \ref{subsection:sparsity}.
	With the polynomial basis defined earlier, the parameterization of $\psi(x_k, z_k, \cdot)$ is now given by the following:
	\begin{align}\label{eq:parameterization}
		\psi(x_k, z_k, \tau) & = q_k p^{*}_{0}(\tau) + \!  \sum\limits_{j \in \mathcal{D}} z_{k,j} p^{*}_{j}(\tau) + hv_k p^{*}_{\mathrm{v}}(\tau).
	\end{align}
	
	Combining \eqref{eq:acceleration} with \eqref{eq:parameterization}, we obtain the PC method explicitly as follows:
	\begin{align}\label{eq:acceleration-explicit}
		\begin{split}
			&x_{k+1} =
			\begin{pmatrix}
				\frac{1}{h} \! \left(q_k \dot{p}^{*}_{0}(1) + \! \sum\limits_{j \in \mathcal{D}} z_{k,j} \dot{p}^{*}_{j}(1) \right) + v_k \dot{p}^{*}_{\mathrm{v}}(1) \\
				q_k p^{*}_{0}(1) + \!  \sum\limits_{j \in \mathcal{D}} z_{k,j} p_{j}(1) + hv_k p^{*}_{\mathrm{v}}(1)
			\end{pmatrix}
			,  \\
			&q_k \ddot{p}^{*}_{0}(\tau_i) \!+\! \! \sum\limits_{j \in \mathcal{D}} \! z_{k, j}\ddot{p}^{*}_{j}(\tau_i) + hv_k \ddot{p}^{*}_{\mathrm{v}}(\tau_i) = \! h^2 f\left(t_k+h\tau_{i}, z_{k,i}\right) \\
			&\hspace{6cm}  \mbox{for all} \;i \in \mathcal{D}.
		\end{split}	
	\end{align}
	
	Note that when $f$ does take $\dot{q}$ as an input, the right-hand side of the last equality has to be replaced with $h^2 f\left(t_k+h\tau_{i}, z_{k,i}, w_{k,i}\right)$ where $w_{k,i}$ are the velocities evaluated on the collocation roots, defined as follows:
	\begin{align}\label{eq:velocities}
		w_{k,i} & \coloneqq  \frac{1}{h} \! \left(q_k \dot{p}^{*}_{0}(\tau_{i}) + \! \sum\limits_{j \in \mathcal{D}} z_{k,j} \dot{p}^{*}_{j}(\tau_{i}) \right) + v_k \dot{p}^{*}_{\mathrm{v}}(\tau_{i}).
	\end{align}
	It is quite clear from that expression that the efficiency improvement brought from the PC methods only takes place when the dependency of $f(\cdot)$ in $\dot{q}$ is small (in the sense of complexity), or nil.

	Let us note that equations \eqref{eq:standard-lagrange} and \eqref{eq:acceleration-explicit} define the functions $F(\cdot)$ and $G(\cdot)$ from the condensed form \eqref{eq:collocation-stacked}.
	A comparison of sparsity and complexity between equations \eqref{eq:standard-lagrange} and \eqref{eq:acceleration-explicit} is provided in Subsection \ref{subsection:sparsity}.

	\section{Direct Collocation for Continuous Optimal Control}\label{section:direct-collocation}
	In this section, we present the class of problems that we aim to tackle and give a description of direct collocation methods which solve them numerically using the collocation methods defined in the previous section.
	
	Now that control is involved, we use the ODE in the form of \eqref{2ndorderODE-control} again, with the same change of variable $x = (\dot{q}, q)$ to end up into an ODE of the following form:
	\begin{align}\label{ODE-control}
		\dot{x}(t) & = \bar{f}(x(t), u(t), p), \quad t \in [0, T].
	\end{align}

	\subsection{Optimal Control Problems}
	We aim to tackle OCP of the form:
	\begin{align}\label{opti:OCP}
		\begin{split}
			& \underset{\substack{x(\cdot) \in \Xx, \\ u(\cdot) \in \Uu, p \in \Pp}}{\mathrm{minimize}} \hspace{0.4cm}
			\int\limits_{0}^{T}  \Ll\left(x(t), u(t), p \right) \dt + E(q(T))  \\
			&\mathrm{subject} \; \mathrm{ to} \\
			& \dot{x}(t) - \bar{f} \!\left( x(t), u(t), p \right) = 0,
			\hspace{1cm} \mbox{for all} \, t \in [0, T], \\
			& g \left(x(t), u(t), p \right) \geq 0,
			\hspace{2cm} \mbox{for all} \, t \in [0, T], \\
			& r(x(T), x(0), p) = 0.
		\end{split}	
	\end{align}

	The spaces $\Xx, \Uu, \Pp$ are respectively the state, control, and parameter spaces.
These are respectively subsets of the set of functions from $[0, T]$ to $\rr^{2\nq}$,
	the set of functions from $[0, T]$ to $\rr^{\nuu}$,
	and $\rr^{\np}$.
	
	The functions $g : \rr^{2\nq} \times \rr^{\nuu} \to \rr^{\nc}$, $\Ll : \rr^{2\nq} \times \rr^{\nuu} \to \rr$, $r :\rr^{2\nq} \times  \rr^{2\nq} \times \rr^{\np}  \to \rr^{\nr} $, and $E : \rr^{\nq} \to \rr$ define the remainder of the control problem.
	The integers $\nq, \nuu, \np, \nc$ and $\nr$ are respectively the dimensions of the state, control, parameter variables, and constraints.
	
	\subsection{Direct collocation methods}
	Direct control methods transform the infinite dimensional OCP \eqref{opti:OCP} into a finite dimension NLP: the controls $u(\cdot)$ are restricted to be piece-wise constant on the equidistant grid defined in the previous section by $N$, $h$, and the grid-points $t_k$ for  $k \in \mathcal{N}$.
	A side remark is that here we consider the case of the same grid for the discretization of $u(\cdot)$ and $x(\cdot)$ which does not have to be the case in general.

	One subclass is the direct transcription methods, which are fully simultaneous approaches where the numerical simulation equations are formulated as equality constraints of the optimization problem, keeping all intermediate variables as optimization variables.
	Probably the most popular class of direct transcription methods is formed by the direct collocation methods, where the integration is performed by collocation method.
	We keep as optimization variables the intermediate variables $z_k$ used in the numerical simulation equations.
	A more detailed description of direct collocation can be found in \cite[Section~8.5.3]{rawlings2017model}.
	An important characteristic of these methods is that it yields a large-scale but sparse  NLP that one can solve with specialized algorithms such as the sparsity exploiting interior-point implementation \textit{IPOPT} \cite{wachter2006implementation}.

	Regarding the implicit integration methods described in the previous section, now that control is involved, they take the following form:
	\begin{align}\label{eq:collocation-stacked-control}
		\begin{split}
			x_{k+1} - F(x_k, z_k) &= 0,  \\
			G(x_k, z_k, u_k, p) &= 0,
		\end{split}	
	\end{align}
	where, consistently with the previous section, $x_k \coloneqq (v_k, q_k) \in \rr^{2\nq}$ approximates $x(t_k) \coloneqq \left( \dot{q}(t_k), q(t_k) \right)$ and $z_k \in \rr^{\nint}$ are the intermediate variables depending on the chosen collocation method and $u_k$ are the values of the controls assumed to be constant on the intervals $[t_k, t_{k+1}]$.
	
	Finally, we use quadrature rules for the integral in the objective of the OCP \eqref{opti:OCP} summarized in the function $l : \rr^{2\nq} \times \rr^{\nint} \to \rr$.
	
	The NLP to be solved in direct collocation for the OCP \eqref{opti:OCP} has as optimization variables:
	the sequences of states $\bm{x} \coloneqq (x_{0}, x_{1}, ..., x_{N})$,
	the sequence of intermediate variables of the collocation equations $\bm{z} \coloneqq \left(z_{0}, ..., z_{N-1}\right)$,
	as well as the sequence of local control parameters $\bm{u} \coloneqq (u_{0}, ..., u_{N-1})$ and $p \in \rr^{ \np}$ 
	and is formulated as follows:
	\begin{align}\label{opti:direct}
		\begin{split}
			& \underset{\bm{x}, \bm{z}, \bm{u}, p}{\mathrm{minimize}} \hspace{0.4cm}
			\sum\limits_{ k \in \mathcal{N}}l\left(x_{k}, z_{k}, u_{k}, p \right) + E(x_{N}) 
			\\
			&\mathrm{subject} \; \mathrm{ to} \\
			& x_{k+1} - F(x_{k}, z_{k}) = 0,
			\hspace{1cm}\mbox{for all} \,k \in \mathcal{N}, \\
			& G(x_{k}, z_{k}, u_{ k}, p)\! = \! 0, \,
			\hspace{1.35cm} \mbox{for all} \, k \in \mathcal{N}, \\
			& g\left( x_{k}, u_{k} \right) \geq 0,
			\hspace{2.15cm} \mbox{for all} \, k \in \mathcal{N}, \\
			& r(x_N, x_{0}, p) = 0.	
		\end{split}
	\end{align}

	Finally, based on these definitions and the previous section,
	we can now make a distinction between the two types of direct methods: the ones with standard collocation for the discretization, and the ones with position-based discretizations.
	The PC methods are obtained by only parameterizing the state variables $q$ as piece-wise polynomial functions whereas SC methods parameterize both the positions and velocities as piece-wise polynomial functions.
	Finally, we recall that in SC methods we have $z_k \in \rr^{2d\cdot \nq}$ because polynomials for both $\dot{q}(\cdot)$ and $q(\cdot)$ are parameterized whereas in PC methods we only parameterize polynomials for $q(\cdot)$, hence $z_k \in \rr^{d\cdot\nq}$ has only half the number of components.
	
	\section{Theoretical comparison}\label{section:theoretical}
	In this Section, we discuss some advantages of using PC methods compared to the standard ones.
	More precisely, we show that, while they give the same order of accuracy, the new approach decreases the complexity of the NLP to solve.
	This provides good reasons to think that it leads to better performance.
	
	\subsection{Order of accuracy}\label{subsection:accuracy}
	
	As it has been shown in \cite{van1991stability} for simulation, the standard and position-based collocation methods, defined in \eqref{eq:standard} and \eqref{eq:acceleration}, converge with the same order.
	More precisely, under sufficient regularity assumptions on the function $f$ in the IVP given by \eqref{2ndorderODE}, the order of convergence of the global error is between $d$ and $2d$ depending on the choice of collocation roots as proven in \cite{van1991stability}.
	Furthermore, it has also been proven in \cite{van1991stability} that in PC methods, the local error on the state $q(\cdot)$ converges with a speed of one order more than in standard methods.
In general, this is however not the case for the error of $\dot{q}(\cdot)$, hence the global error has a convergence of the same order as the standard collocation and not more.
	
	This is however only due to the approximation error on $\dot{q}(\cdot)$, which gives us reasons to expect a convergence with a better multiplicative constant, even though the order is the same.
	
\begin{figure}
		\begin{center}
			\scalebox{0.65}{
				\input{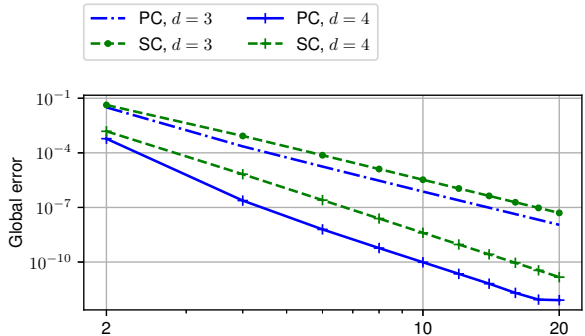}
			}
		\end{center}
		\vspace{-0.8cm}
		\caption{Accuracy of the two methods on the IVP \eqref{eq:example} for Legendre collocation points.
			A collocation stands for PC method and S collocation stands for standard collocation method.
			The global error is the $L_1$ distance to the analytical solution.
		}
		\label{fig:collocation-accuracy}
	\end{figure}

	As an alternative interpretation, the polynomial which locally approximates the solution $q(\cdot)$ is one degree higher than in the PC method \eqref{eq:acceleration}.
	As an example, the PC method of order $d$ with collocation roots of any choice would offer an exact solution for the ODE $\ddot{q}(t) = t^{d-1}$, while the standard collocation method of order $d$ would not.
	The convergence speed of the error is illustrated in Figure \ref{fig:collocation-accuracy} for the following example IVP:
		\begin{align}
			\label{eq:example}
			\begin{split}
				\ddot{q}(t) + q(t) &=  \cos(t), \; \mbox{for all} \; t\in [0, 10], \\
				q(0) = \dot{q}(0) & = 0.
			\end{split}
		\end{align}
	The plots show that for this example, the PC method converges with the same order (the order is $2d$ here because of the choice of collocation roots of Gauss-Legendre type), but with a better multiplicative constant.
	
	\subsection{Number of variables and sparsity pattern}\label{subsection:sparsity}
	
	In this section, we assume that the functions $L(\cdot)$ and $f(\cdot)$ in \eqref{opti:OCP} only depend on $q$ and not on $\dot{q}$.
	This restriction is not necessary for the PC methods to benefit from the same order of convergence, but complexity reduction is only interesting in that case.
Let us also note that for given $N$ and $d$, the number of evaluations of functions $L(\cdot)$, $E(\cdot)$, $f(\cdot)$, $g(\cdot)$, and $r(\cdot)$ is the same for the two methods.
	Finally, in the two methods, the number of evaluations of the latter functions is the same.
	
	Let us study the dimension and sparsity of the NLP obtained from PC and SC methods.
	Note that the number of nonlinear function evaluations of the constraints and objective and the number of structural nonzero elements in the gradient of the objective function are the same for the two methods.
	
	In Table \ref{table:complexity}, we compare the number of variables, constraints, and nonzero elements in the constraint-Jacobian, for the PC and SC methods.
	In this table, we used the following constants, for the values that are independent of the method:
	\begin{align*}
		\begin{split}
			C_{1} &\coloneqq N \left(2\nq + \nuu \right) +\np + 2\nq
			, \\
			C_2 &\coloneqq 2\nq + \nr + N \nc + 2N \nq 
			,\\
			C_3 &\coloneqq N\nq d(\nq + \nuu + \np)+ N\nc (\nuu + 2\nq) \nonumber \\
			& \hspace{4cm}  + 2\nq (2\nq + \nr).	
		\end{split}
	\end{align*}
	The number $C_1$ is for the number of variables, $C_2$ is for the number of constraints and $C_3$ is for the number of nonzero elements in constraints-Jacobian.
	
	It happens that the PC method has a sparser and smaller constraint-Jacobian.
	This simply comes from the fact that in SC methods, there are twice as many intermediate variables and yet, no special sparsity is gained when $L(\cdot)$ and $f(\cdot)$ are independent of $\dot{q}$.

	\begin{table}
		\caption{Complexity comparison of direct collocation methods}
		\begin{tabular}{ |L{2.3cm}|L{2.45cm}|L{2.45cm}| }
			\hline && \\
			& \textbf{SC method} & \textbf{PC method}
			\\ \hline && \\
			\textbf{Number of variables} & $C_1 + 2 N \nq d $ &  $C_1 + N \nq d$
			\\ \hline && \\
			\textbf{Number of constraints} & $C_2 + 2N \nq d $ & $C_2 + N \nq  d $
			\\ \hline && \\
			\textbf{Nonzero elements in constraint-Jacobian} & $C_3$    $+N\nq(2d^2+5d+4)$ & $C_3$   $+N\nq(d^2+3d+6)$ 
			\\ \hline
		\end{tabular}
	\label{table:complexity}
	\end{table}

	\section{Numerical comparison}\label{section:numerical}
	In this section, numerical results concerning the described method are shown.
These are based on a popular example of a second-order OCP.
	
	\begin{figure}[h!]
		\begin{center}
			\scalebox{0.8}{
				\input{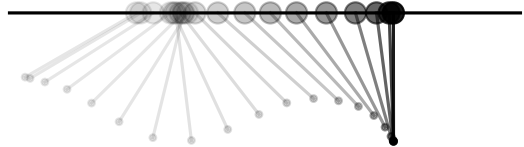}
			}
		\end{center}
		\caption{Trajectory of a two-dimensional overhead crane system.}
		\label{fig:crane}
	\end{figure}
	
	\subsection{Model description}
	
	The example we choose for comparison of PC and SC methods is a problem of optimal control of a frictionless overhead crane with inextensible rope.
	The equations of such a system are described in \cite{wu2015model} for example.
	However, in our case, we assume the mass of the load to be infinitely smaller than the mass of the trolley, such that the motion of the trolley is only influenced by the control force, and not by the tension of the rope.
	We only model a fluid-friction torque applied to the rotation of the rope.
	Equivalently, it can be seen as directly controlling the acceleration of the trolley instead of the force, which also justifies the absence of friction for the trolley in the model, which simplifies the equations.
	We perform optimal control with box constraints on positions and terminal constraints on both positions and velocities for the equilibrium point at the origin.
	We minimize the squared norm of the controls and positions along the trajectory.
	This problem is also called a regulation problem.
	The state is $q(t) = \left( r(t), \theta(t) \right)$ where $r(t)$ represents the position of the trolley, and $\theta(t)$ represents the angle of the rope at time $t$.
	A controlled trajectory for this problem is illustrated in Figure \ref{fig:crane}, and the corresponding OCP is the following:
	
	\begin{align}\label{opti:ex}
		\begin{split}
			& \underset{\substack{r(\cdot), \theta(\cdot) ,u(\cdot)}}{\mathrm{minimize}} \hspace{0.4cm}
			\int_{0}^{T} \left(u(t)^2 +  r(t)^2 + \theta(t)^2  \right)  \,  \dt
			\\
			&\mathrm{subject} \; \mathrm{ to} \\
			& \ddot{r}(t) =  \,  u(t),
			\hspace{3.1cm} \mbox{for all} \;  t \in [0, T], \\
			& \ddot{\theta}(t) = f\left(\theta(t), u(t)   \right) - \beta \dot{\theta},
			\hspace{1.1cm} \mbox{for all} \;  t \in [0, T], \\
			& r_{\min} \leq r(t)  \leq r_{\max},
			\hspace{2cm} \mbox{for all} \;  t \in [0, T], \\
			& \big( \dot{r}(0), \dot{\theta}(0)  \big) = \left(0, 0\right), \\
			& \big( r(0), \theta(0)  \big) = \left(r_0, \theta_0\right), \\
			& \big( \dot{r}(T), \dot{\theta}(T)  \big) = \left(0, 0\right), \\
			& \big( r(T), \theta(T)  \big) = \left(0, 0\right),
		\end{split}
	\end{align}
	where the function $f(\cdot)$ is defined as follows:
	\begin{align}\label{eq:ex}
		f\left(\theta, u\right) &\coloneqq - u \cos \left(\theta \right)  - a \sin \left(\theta\right).
	\end{align}
	
	
	To illustrate the adaptability of the method, we use both Gauss-Legendre and Radau IIA collocation points \cite{rawlings2017model} for both methods.
	We implement the methods described in this paper in \textit{Python}, using the package  \textit{CasADi} \cite{Andersson2019} for manipulation of symbolic variables.
	To solve the NLP, we use \textit{IPOPT} \cite{wachter2006implementation} via its \textit{CasADi} interface in \textit{Python}.
	The results are presented in the next part.

	\subsection{Numerical results}
	To compare the accuracy of both methods, we apply them on $200$ instances of the problem \eqref{opti:ex}, which are generated by drawing $(r_0, \theta_0)$ from a uniform distribution on $[-3, 3] \times [-\frac{\pi}{3}, \frac{\pi}{3}]$.
	For the rest of the parameters of \eqref{opti:ex} and equation \eqref{eq:ex}, we use $r_{\min} = -3.0$, $r_{\max} = 3.0$, $T = 10.0$, $\beta=0.1$ and $a = 9.81$.
	Furthermore, we apply the PC and SC methods with $N = 20$ and for two different values of the order: $d = 2$ and $d = 3$.
	To compare performance, for a given instance of the problem and a given method, we compute the error as follows: $\mathrm{error} \coloneqq v - v^{\mathrm{ref}} $ where $v$ is the objective value of the problem evaluated for the solution given by the method of interest, while $v^{\mathrm{ref}}$ is the objective value when the solution is computed with higher accuracy, using SC with $d = 5$.
	On the other hand, to compare the efficiency of the methods, we compare the average running time per iteration of the interior-point method performed through \textit{IPOPT}.
	This approach is preferred to total computational time because it ignores the number of iteration which depends on the initial guess.
	The code is executed on Python 3.8.8 on a standard laptop.
	
	\begin{figure}
		\begin{center}
			\scalebox{0.8}{
				\input{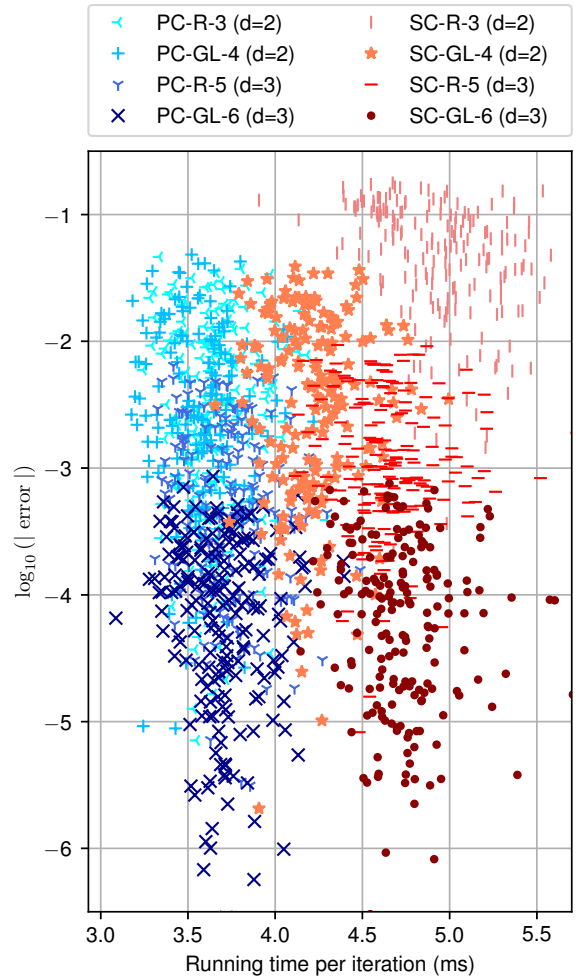}
			}
		\end{center}
		\vspace{-0.5cm}
		\caption{Pareto plot:
			the log of the error for each problem is plotted against the average running time per iteration.
		}
		\label{fig:pareto}
	\end{figure}
	\begin{table}
		\vspace{0.5cm}
		\caption{Comparison of the two methods for the friction-less example}
		\begin{tabular}{|L{0.6cm}|C{1.2cm}|C{2cm}|C{2cm}|}
			\hline
			& & & 
			\\
			& & \textbf{Average running time per iterations in milliseconds} & \textbf{Geometric mean of the error}
			\\ \hline & & &  \\
			& \textbf{$d=2$ Radau} & $5.0$ &  $4.5 \cdot 10^{-2}$
			\\ \cline{2-4} & & &  \\
			\textbf{SC}  & \textbf{$d=2$ Legendre} & $4.2$ &  $3.0 \cdot 10^{-3}$
			\\ \cline{2-4} & & &  \\
			& \textbf{$d=3$ Radau} & $4.7$ & $1.2 \cdot 10^{-3}$
			\\ \cline{2-4} & & &  \\
			& \textbf{$d=3$ Legendre} & $4.7$ & $5.6 \cdot 10^{-5}$
			\\ \hline & & &  \\
			& \textbf{$d=2$ Radau} & $3.7$ &  $2.5 \cdot 10^{-3}$
			\\ \cline{2-4} & & &  \\
			\textbf{PC} & \textbf{$d=2$ Legendre} & $3.6$ &  $2.7 \cdot 10^{-3}$
			\\ \cline{2-4} & & &  \\
			& \textbf{$d=3$ Radau} &$3.7$ &  $4.6 \cdot 10^{-4}$
			\\ \cline{2-4} & & &   \\
			& \textbf{$d=3$ Legendre} &$3.6$ &  $6.3 \cdot 10^{-5}$
			\\ \hline
		\end{tabular}
	\label{table:results}
	\end{table}
	
	Figure \ref{fig:pareto} and Table \ref{table:results}, summarize the result.
	We use the notation "X-Y-$n$" for the method X with the collocation roots given by Y, which has a convergence of order $n$.
	We use the abbreviation GL for Gauss-Legendre collocation points and R for Radau-IIA.
	We observe that, in this example, the PC and SC methods give solutions with similar quality, but the PC methods are faster than SC methods, by $15$ to $25 \%$ in running time.
	We observe this gain of efficiency both in the cheap methods ($d=2$) and in the accurate methods ($d=3$), and for both choices of collocation points.
	This example show that even for an ODE where the right-hand side has a dependency in the velocities, the PC methods are more efficient.
	
	\section{Conclusion}
	This paper presents and compares different variations of direct collocation methods for second-order ODE.
With these comparisons, we have shown that in the presented simple example, the PC method seems to present an advantage in terms of running time compared to the standard ones, but not significantly in terms of accuracy.
	Hence, we can conclude that our theoretical results from Section \ref{section:theoretical} have concrete implications for the application of our method to real-world problems.
	Indeed, even if the generic gain in sparsity is obtained when the right-hand side of the ODE is independent of the velocities, the advantage is still visible in the presented example, even though the right-hand side does depend on the velocities there.
	This can be interpreted by the fact that this dependency is light (linear, only on some part of the components).
	
	Finally, these methods should also be used and studied for problems of the same kind with higher dimensions, typically for applications from robotics.

	\bibliographystyle{IEEEtran}  
	\bibliography{bibli}

\begin{thebibliography}{10}
\providecommand{\url}[1]{#1}
\csname url@samestyle\endcsname
\providecommand{\newblock}{\relax}
\providecommand{\bibinfo}[2]{#2}
\providecommand{\BIBentrySTDinterwordspacing}{\spaceskip=0pt\relax}
\providecommand{\BIBentryALTinterwordstretchfactor}{4}
\providecommand{\BIBentryALTinterwordspacing}{\spaceskip=\fontdimen2\font plus
\BIBentryALTinterwordstretchfactor\fontdimen3\font minus
  \fontdimen4\font\relax}
\providecommand{\BIBforeignlanguage}[2]{{%
\expandafter\ifx\csname l@#1\endcsname\relax
\typeout{** WARNING: IEEEtran.bst: No hyphenation pattern has been}%
\typeout{** loaded for the language `#1'. Using the pattern for}%
\typeout{** the default language instead.}%
\else
\language=\csname l@#1\endcsname
\fi
#2}}
\providecommand{\BIBdecl}{\relax}
\BIBdecl

\bibitem{magni2009nonlinear}
L.~Magni, D.~M. Raimondo, and F.~Allg{\"o}wer, ``Nonlinear model predictive
  control,'' \emph{Lecture Notes in Control and Information Sciences}, vol.
  384, 2009.

\bibitem{rawlings2017model}
J.~B. Rawlings, D.~Q. Mayne, and M.~Diehl, \emph{Model predictive control:
  theory, computation, and design}.\hskip 1em plus 0.5em minus 0.4em\relax Nob
  Hill Publishing Madison, 2017, vol.~2.

\bibitem{lurie2002analytical}
A.~I. Lurie, \emph{Analytical mechanics}.\hskip 1em plus 0.5em minus
  0.4em\relax Springer Science \& Business Media, 2002.

\bibitem{winkler2018gait}
A.~W. Winkler, C.~D. Bellicoso, M.~Hutter, and J.~Buchli, ``Gait and trajectory
  optimization for legged systems through phase-based end-effector
  parameterization,'' \emph{IEEE Robotics and Automation Letters}, vol.~3,
  no.~3, pp. 1560--1567, 2018.

\bibitem{kramarz1980stability}
L.~Kramarz, ``Stability of collocation methods for the numerical solution of
  y''= f (x, y),'' \emph{BIT Numerical Mathematics}, vol.~20, no.~2, pp.
  215--222, 1980.

\bibitem{van1991stability}
P.~Van~der Houwen, B.~Sommeijer, and N.~H. Cong, ``Stability of
  collocation-based runge-kutta-nystr{\"o}m methods,'' \emph{BIT Numerical
  Mathematics}, vol.~31, no.~3, pp. 469--481, 1991.

\bibitem{d2009two}
R.~D’Ambrosio, M.~Ferro, and B.~Paternoster, ``Two-step hybrid collocation
  methods for y''= f (x, y),'' \emph{Applied Mathematics Letters}, vol.~22,
  no.~7, pp. 1076--1080, 2009.

\bibitem{martucci2005general}
S.~Martucci and P.~Beatrice, ``General two step collocation methods for special
  second order ordinary differential equations,'' in \emph{Paper Proceedings of
  the 17th IMACS World Congress Scientific Computation, Applied Mathematics and
  Simulation, Paris}, 01 2005.

\bibitem{moreno2022collocation}
C.~Moreno-Mart\'in, Ros, ``Collocation methods for second order systems,''
  \emph{Robotics: Science and Systems}, vol.~18, 2022.

\bibitem{wachter2006implementation}
A.~W{\"a}chter and L.~T. Biegler, ``On the implementation of an interior-point
  filter line-search algorithm for large-scale nonlinear programming,''
  \emph{Mathematical programming}, vol. 106, no.~1, pp. 25--57, 2006.

\bibitem{wu2015model}
Z.~Wu, X.~Xia, and B.~Zhu, ``Model predictive control for improving operational
  efficiency of overhead cranes,'' \emph{Nonlinear Dynamics}, vol.~79, no.~4,
  pp. 2639--2657, 2015.

\bibitem{Andersson2019}
J.~A.~E. Andersson, J.~Gillis, G.~Horn, J.~B. Rawlings, and M.~Diehl,
  ``{CasADi} -- {A} software framework for nonlinear optimization and optimal
  control,'' \emph{Mathematical Programming Computation}, vol.~11, no.~1, pp.
  1--36, 2019.

\end{thebibliography}
	
\end{document}